\documentclass[12pt,draft,leqno]{article}
\usepackage{amssymb,eucal,latexsym}
\usepackage{amsmath,amsthm,latexsym,amscd,amsbsy,fleqn,graphicx}

\newtheorem{theorem}{Theorem}[section]
\newtheorem{lm}[theorem]{Lemma}
\newtheorem{exa}[theorem]{Example}

\newtheorem{pro}[theorem]{Proposition}
\newtheorem{defi}[theorem]{Definition}

\newtheorem{nota}[theorem]{Notation}

\newtheorem{rem}[theorem]{Remark}

\newtheorem{fact}[theorem]{Fact}

\def\p{\varphi}
\def\a{\alpha}

\def\d{\delta}

\def\g{\gamma}
\def\GA{\Gamma}
\def\DE{\Delta}

\def\l{\lambda}

\def\s{\sigma}

\def\lag{\lambda_A^g}
\def\lbg{\lambda_B^g}

\def\lra{\longrightarrow}

\def\sbe{\subseteq}

\def\stm{\setminus}
\def\ems{\emptyset}
\def\nes{\neq\emptyset}

\def\wt{\widetilde}

\def\unl{\underline}

\def\ex{\exists}
\def\fa{\forall}

\def\we{\wedge}

\def\bv{\bigvee}

\def\ap{^{\,\prime}}
\def\inv{^{-1}}
\def\st{\ |\ }

\def\nin{\not\in}

\def\card #1{\vert #1 \vert}

\def\CC{{\cal C}}

\def\TT{{\cal T}}

\def\HLC{{\bf HLC}}

\def\MDHLC{{\bf MDHLC}}

\def\1{{\bf 1}}
\def\2{\mbox{{\bf 2}}}
\def\3{\mbox{{\bf 3}}}

\def\int{\mbox{{\rm int}}}

\def\cl{\mbox{{\rm cl}}}
\def\CL{\mbox{{\rm Clust}}}
\def\Clust{\mbox{{\rm Clust}}}
\def\BClu{\mbox{{\rm BClust}}}

\def\doc{\hspace{-1cm}{\em Proof.}~~}
\def\sq{\hspace*{\fill} \hbox{\vrule\vbox{\hrule\phantom{o}\hrule}\vrule}}
\def\sqs{\sq \vspace{2mm}}

\def\BBBB{\mathbb{B}}

\def\NNNN{\mathbb{N}}

\def\di{\diamond}

\def\bU0{\bar{U}=(U^0,(U^i,U^{ci})_{i\in\omega})}

\def\bV0{\bar{V}=(V^0,(V^i,V^{ci})_{i\in\omega})}

\title{{\huge\bf A new duality theorem for locally}\\
\vspace{0.2cm}
{\huge\bf compact spaces}\\
\vspace{0.5cm}
{\large\bf Georgi Dimov, Elza Ivanova}\\
\vspace{0.2cm}
{\footnotesize\rm Department of Mathematics and Informatics, Sofia University, 5 J. Bourchier Blvd., 1164 Sofia,}\\
{\footnotesize\rm Bulgaria}
}

\author{}

\date{}

\begin{document}

\maketitle

\begin{abstract}
{\footnotesize
\noindent In 1962, de Vries \cite{deV} proved a duality theorem for the
category {\bf HC} of compact Hausdorff spaces and continuous maps. The composition of the morphisms of the dual category obtained by him
differs from the set-theoretic one.
 Here we obtain a new category dual to the category {\bf HLC} of locally compact Hausdorff spaces and
continuous maps for which the composition of the morphisms is a natural one but the morphisms are multi-valued maps.}
\end{abstract}

{\footnotesize {\em  MSC:} primary 54D45, 18A40; secondary 54E05.

{\em Keywords:} Locally compact space; Local contact algebra;
   Continuous map; $\delta$-ideal; Duality.}

\footnotetext[1]{{\footnotesize{This paper was supported by the
project no. 173/2010 $``$Duality, equivalence and representation theorems and their applications in general topology" of the Sofia
University $``$St. Kl. Ohridski".}}}

\footnotetext[2]{{\footnotesize {\em E-mail addresses:}
gdimov@fmi.uni-sofia.bg, elza@fmi.uni-sofia.bg}}

%-------------------------------------------------------------------------
%-----------------------------Introduction--------------------------------
%--------------------------------------------------------------------------

\section{Introduction}

In 1962, de Vries \cite{deV} proved a duality theorem for the
category {\bf HC} of compact Hausdorff spaces and continuous maps.
This theorem was the first realization in a full extent of the
ideas of the so-called {\em region-based theory of space},
although, as it seems, de Vries did not know of the existence of
such a theory. The region-based theory of space is a kind of
point-free geometry and can be considered as  an alternative to
the well known Euclidean point-based theory of space. Its main
idea goes back to Whitehead \cite{W} (see also \cite{W1}) and de Laguna \cite{dL}
and is based on a certain criticism of the Euclidean approach to
the geometry, where the points (as well as straight lines and
planes) are taken as the basic primitive notions. A. N. Whitehead and T. de Laguna noticed that
points, lines and planes are quite abstract entities which have
not a separate existence in reality and proposed to put the theory
of space  on the base of some more realistic spatial entities. In
Whitehead \cite{W}, the notion of region is taken as a primitive
notion: it is an abstract analog of a spatial body; also some
natural relations between regions are regarded. In \cite{W1} Whitehead
considers only some mereological relations like ``part-of" and
``overlap", while in \cite{W} he adopts  from de Laguna \cite{dL} the
relation of ``contact" (``connectedness" in Whitehead's original
terminology) as the only primitive relation between regions. In
this way the region-based theory of space appeared as an extension
of mereology -- a philosophical discipline of ``parts and wholes".

 Let us note that neither A. N. Whitehead nor T. de
Laguna presented  their ideas in a detailed mathematical form. Their
ideas attracted some mathematicians and mathematically oriented
philosophers to present various versions of region-based theory of
space at different levels of abstraction. Here we can mention
A. Tarski \cite{Tarski}, who rebuilt Euclidean geometry as an extension of
mereology with the primitive notion of sphere. Remarkable is also
Grzegorczyk's paper \cite{Grzegorczyk}. Models of Grzegorczyk's  theory  are complete
Boolean algebras of regular closed sets of certain topological
spaces equipped with the relation of {\em separation}\/ which in fact is
the complement of Whitehead's contact relation. On the same line
of abstraction is also the point-free topology \cite{J}.

Let us mention that Whitehead's ideas of region-based theory of
space  flourished and in a sense  were reinvented and applied in
some areas of computer science: Qualitative Spatial Reasoning
(QSR), knowledge representation, geographical information systems,
formal ontologies in information systems, image processing,
natural language semantics etc. The reason is that the language of
region-based theory of space allows us to obtain a more simple
description of some qualitative spatial features and properties of
space bodies.  One of the most popular among the community of
QSR-researchers is the system of Region Connection Calculus (RCC)
introduced by Randell, Cui and Cohn \cite{Randell}.

A celebrated
duality for the category {\bf HC} is the Gelfand Duality Theorem \cite{G1,G2,GN, GS}.
 The de Vries Duality Theorem, however, is
the first  complete realization of the  ideas of de Laguna \cite{dL}
and Whitehead \cite{W}: the models of the regions in de Vries' theory
are the regular closed sets of compact Hausdorff spaces (regarded with the well known  Boolean
structure on them) and the contact relation $\rho$ between these sets is defined by $F\rho G \iff F\cap G\nes$.

The composition of the morphisms of  de Vries' category {\bf DHC}
dual to the category {\bf HC} differs from their set-theoretic
composition. In 1973, V. V. Fedorchuk \cite{F} noted that the complete
{\bf DHC}-morphisms (i.e., those {\bf DHC}-morphisms which are
complete Boolean homomorphisms) have a very simple description
and, moreover, the {\bf DHC}-com\-po\-sition of two such morphisms
coincides with their set-theoretic composition. He considered the
{\em cofull subcategory} (i.e. such a subcategory which has the
same objects as the whole category) {\bf DQHC} of the category
{\bf DHC} determined by the complete {\bf DHC}-morphisms. He
proved that the restriction of de Vries' duality functor to it
produces a duality between the category {\bf DQHC} and the
category {\bf QHC} of compact Hausdorff spaces and quasi-open maps
(a class of maps introduced by Marde\v{s}ic and Papic in \cite{MP}).

It is natural to try to extend de Vries' Duality Theorem to the
category {\bf HLC} of locally compact Hausdorff spaces and
continuous maps. An important step in this direction was done by
Roeper \cite{R}. Being guided by the ideas of de Laguna \cite{dL} and
Whitehead \cite{W}, he defined the notion of {\em region-based
topology} which is now known as {\em local contact algebra}\/
(briefly, LCA or LC-algebra) (see \cite{DV}), because the axioms which it
satisfies almost coincide with the axioms of local proxi\-mities
of Leader \cite{LE}. In his paper \cite{R}, Roeper proved the following
theorem: there is a bijective correspondence between all (up to
homeomorphism) locally compact Hausdorff spaces and all (up to
isomorphism) complete LC-algebras. In \cite{D2009}, using Roeper's theorem, the Fedorchuk Duality Theorem
was extended to the category of locally compact Hausdorff spaces and skeletal (in the sense of \cite{MR}) maps.
Quite recently,
in the paper \cite{D-AMH1-10},  de Vries' Duality Theorem \cite{deV} was
extended to the category {\bf HLC}. The composition of the morphisms of
the obtained there dual category is not the usual composition of
maps (i.e., the situaton is the same as in the case of de Vries' Duality Theorem). We now
obtain a new duality theorem for the category {\bf HLC} such that
the composition of the morphisms of the dual category is a natural
one (like in the Fedorchuk Duality Theorem for the category {\bf QHC}); however, the morphisms of the dual category are multi-valued
maps.

Let us fix the notation.

If $\CC$ denotes a category, we write $X\in \card\CC$ if $X$ is an
object of $\CC$, and $f\in \CC(X,Y)$ if $f$ is a morphism of $\CC$
with domain $X$ and codomain $Y$.

All lattices are with top (= unit) and bottom (= zero) elements,
denoted respectively by 1 and 0. We do not require the elements
$0$ and $1$ to be distinct. The operation $``$complement" in
Boolean algebras is denoted by $``$*". The (positive) natural numbers are
denoted by $\mathbb{N}$ (resp., by $\mathbb{N}^+$). The Alexandroff (one-point) compactification of a locally compact Hausdorff space
$X$ is denoted by $\a X$. If $X$ is a set then we will denote by $id_X$ the identity function on $X$.

%---------------------------------------------------------------------------------
\section{Preliminaries}
%--------------------------------------------------------------------------------

\begin{defi}\label{conalg}
\rm
An algebraic system $(B,0,1,\vee,\wedge, {}^*, C)$ is called a {\it
contact Boolean algebra}\/ or, briefly, {\it contact algebra}
(abbreviated as CA or C-algebra) (\cite{DV})
 if the system
$(B,0,1,\vee,\wedge, {}^*)$ is a Boolean algebra
  and $C$
is a binary relation on $B$, satisfying the following axioms:

\smallskip

\noindent (C1) If $a\not= 0$ then $aCa$;\\
(C2) If $aCb$ then $a\not=0$ and $b\not=0$;\\
(C3) $aCb$ implies $bCa$;\\
(C4) $aC(b\vee c)$ iff $aCb$ or $aCc$.

\smallskip

\noindent We shall simply write $(B,C)$ for a contact algebra. The
relation $C$  is called a {\em  contact relation}. When $B$ is a
complete Boolean algebra, we will say that $(B,C)$ is a {\em
complete contact Boolean algebra}\/ or, briefly, {\em complete
contact algebra} (abbreviated as CCA or CC-algebra).

We will say that two C-algebras $(B_1,C_1)$ and $(B_2,C_2)$ are  {\em
CA-isomorphic} iff there exists a Boolean isomorphism $\varphi:B_1\longrightarrow
B_2$ such that, for each $a,b\in B_1$, $aC_1 b$ iff $\varphi(a)C_2
\varphi(b)$. Note that in this paper, by a $``$Boolean isomorphism" we
understand an isomorphism in the category {\bf Bool} of Boolean algebras and Boolean homomorphisms.

\smallskip

A contact algebra $(B,C)$ is called a {\it  normal contact Boolean
algebra}\/ or, briefly, {\it  normal contact algebra} (abbreviated
as NCA or NC-algebra) (\cite{deV},\cite{F}) if it satisfies the following axioms (we
will write $``-C$" for $``not\ C$"):

\smallskip

\noindent (C5) If $a(-C)b$ then $a(-C)c$ and $b(-C)c^*$ for some $c\in B$;\\
(C6) If $a\not= 1$ then there exists $b\not= 0$ such that
$b(-C)a$.

\smallskip

\noindent Note that the axioms of NC-algebras are very similar to the
Efremovi\v c axioms of proximity spaces \cite{EF}.

 A normal CA is called a {\em complete normal contact
Boolean algebra}\/ or, briefly, {\em complete normal contact
algebra} (abbreviated as CNCA or CNC-algebra) if it is a CCA. The notion of
normal contact algebra was introduced by Fedorchuk \cite{F} under
the name {\em Boolean $\delta$-algebra}\/ as an equivalent expression
of the notion of compingent Boolean algebra of de Vries (see its definition below). We call
such algebras $``$normal contact algebras" because they form a
subclass of the class of contact algebras and naturally arise in
normal Hausdorff spaces.

Note that if $0\neq 1$ then the axiom (C2) follows from the axioms
(C6) and (C4).

For any CA $(B,C)$, we define a binary relation  $``\ll_C $"  on
$B$ (called {\em non-tangential inclusion})  by $``\ a \ll_C b
\leftrightarrow a(-C)b^*\ $". Sometimes we will write simply
$``\ll$" instead of $``\ll_C$".
\end{defi}

The relations $C$ and $\ll$ are inter-definable. For example,
normal contact algebras could be equivalently defined (and exactly
in this way they were introduced (under the name of {\em
compingent Boolean algebras}) by de Vries in \cite{deV}) as a pair
of a Boolean algebra $B=(B,0,1,\vee,\wedge,{}^*)$ and a binary
relation $\ll$ on $B$ subject to the following axioms:

\smallskip

\noindent ($\ll$1) $a\ll b$ implies $a\leq b$;\\
($\ll$2) $0\ll 0$;\\
($\ll$3) $a\leq b\ll c\leq t$ implies $a\ll t$;\\
($\ll$4) $a\ll c$ and $b\ll c$ implies $a\vee b\ll c$;\\
($\ll$5) If  $a\ll c$ then $a\ll b\ll c$  for some $b\in B$;\\
($\ll$6) If $a\neq 0$ then there exists $b\neq 0$ such that $b\ll
a$;\\
($\ll$7) $a\ll b$ implies $b^*\ll a^*$.

\smallskip

Note that if $0\neq 1$ then the axiom ($\ll$2) follows from the
axioms ($\ll$3), ($\ll$4), ($\ll$6) and ($\ll$7).

\smallskip

Obviously, contact algebras could be equivalently defined as a
pair of a Boolean algebra $B$ and a binary relation $\ll$ on $B$
subject to the  axioms ($\ll$1)-($\ll$4) and ($\ll$7).

\smallskip

It is easy to see that axiom (C5) (resp., (C6)) can be stated
equivalently in the form of ($\ll$5) (resp., ($\ll$6)).

\begin{exa}\label{rct}
\rm Recall that a subset $F$ of a topological space $(X,\tau)$ is
called {\em regular closed}\/ if $F=\cl(\int (F))$. Clearly, $F$
is regular closed iff it is the closure of an open set.

For any topological space $(X,\tau)$, the collection $RC(X,\tau)$
(we will often write simply $RC(X)$) of all regular closed subsets
of $(X,\tau)$ becomes a complete Boolean algebra
$(RC(X,\tau),0,1,\wedge,\vee,{}^*)$ under the following operations:
$ 1 = X,  0 = \emptyset, F^* = \cl(X\setminus F), F\vee G=F\cup G,
F\wedge G =\cl(\int(F\cap G)).
$
The infinite operations are given by the  formulas:
$\bigvee\{F_{\gamma}\st \gamma\in\Gamma\}=\cl(\bigcup\{F_{\gamma}\st
\gamma\in\Gamma\}),$ and
$\bigwedge\{F_{\gamma}\st \gamma\in\Gamma\}=\cl(\int(\bigcap\{F_{\gamma}\st
\gamma\in\Gamma\})).$

It is easy to see that setting $F \rho_{(X,\tau)} G$ iff $F\cap
G\not = \emptyset$, we define a contact relation $\rho_{(X,\tau)}$ on
$RC(X,\tau)$; it is called a {\em standard contact relation}. So,
$(RC(X,\tau),\rho_{(X,\tau)})$ is a CCA (it is called a {\em
standard contact algebra}). We will often write simply $\rho_X$
instead of $\rho_{(X,\tau)}$. Note that, for $F,G\in RC(X)$,
$F\ll_{\rho_X}G$ iff $F\subseteq \int_X(G)$.

Clearly, if $(X,\tau)$ is a normal Hausdorff space then the
standard contact algebra $(RC(X,\tau),\rho_{(X,\tau)})$ is a
complete NCA.

A subset $U$ of $(X,\tau)$ such that $U=\int(\cl(U))$ is said to
be {\em regular open}.
\end{exa}

\begin{defi}\label{clust}
\rm
Let $(B,C)$ be a $CA$. Then a non-empty subset $\sigma$ of $B$ is called a
cluster in $(B,C)$ if the following conditions are satisfied:

\noindent $(K1)$ If $a,b\in\sigma$ then $aCb$;

\noindent $(K2)$ If $a\vee b\in\sigma$ then $a\in\sigma$ or $b\in\sigma$;

\noindent $(K3)$ If $aCb$ for every $b\in\sigma$, then $a\in\sigma$.

\noindent The set of all clusters in $(B,C)$ will be denoted by $\Clust (B,C)$.
\end{defi}

\begin{pro}\label{Fed1}{\rm (\cite{D2009}, \cite{R})}
Let $(B,C)$ be a normal contact algebra, $\sigma$ be a
cluster in $(B,C)$, $a\in B$ and $a\not\in\sigma$. Then there exists $b\in B$ such that $b\not\in\sigma$ and $a\ll b$.
 \end{pro}

The following notion is a lattice-theoretical counterpart of the
Leader's notion of {\em local proximity} (\cite{LE}):

\begin{defi}\label{locono}{\rm (\cite{R})}
\rm An algebraic system $$\underline {B}_{\, l}=(B,0,1,\vee,\wedge,
{}^*, \rho, \mathbb{B})$$ is called a {\it local contact Boolean
algebra}\/ or, briefly, {\it local contact algebra} (abbreviated
as LCA or LC-algebra)   if $(B,0,1, \vee,\wedge, {}^*)$ is a Boolean algebra,
$\rho$ is a binary relation on $B$ such that $(B,\rho)$ is a CA,
and $\mathbb{B}$ is an ideal (possibly non proper) of $B$, satisfying
the following axioms:

\smallskip

\noindent(BC1) If $a\in\mathbb{B}$, $c\in B$ and $a\ll_\rho c$ then
$a\ll_\rho b\ll_\rho c$ for some $b\in\mathbb{B}$;\\
(BC2) If $a\rho b$ then there exists an element $c$ of $\mathbb{B}$
such that
$a\rho (c\wedge b)$;\\
(BC3) If $a\neq 0$ then there exists  $b\in\mathbb{B}\setminus\{0\}$ such
that $b\ll_\rho a$.

\smallskip

We shall simply write  $(B, \rho,\mathbb{B})$ for a local contact
algebra.  We will say that the elements of $\mathbb{B}$ are {\em
bounded} and the elements of $B\setminus \mathbb{B}$  are  {\em unbounded}.
When $B$ is a complete Boolean algebra,  the LCA $(B,\rho,\mathbb{B})$
is called a {\em complete local contact Boolean algebra}\/ or,
briefly, {\em complete local contact algebra} (abbreviated as
CLCA or CLC-algebra).

We will say that two local contact algebras $(B,\rho,\mathbb{B})$ and
$(B_1,\rho_1,\mathbb{B}_1)$ are  {\em LCA-isomorphic} if there exists a
Boolean isomorphism $\varphi:B\longrightarrow B_1$ such that, for $a,b\in B$,
$a\rho b$ iff $\varphi(a)\rho_1 \varphi(b)$, and $\varphi(a)\in\mathbb{B}_1$ iff
$a\in\mathbb{B}$.
\end{defi}

\begin{rem}\label{conaln}
\rm Note that if $(B,\rho,\mathbb{B})$ is a local contact algebra and
$1\in\mathbb{B}$ then $(B,\rho)$ is a normal contact algebra.
Conversely, any normal contact algebra $(B,C)$ can be regarded as
a local contact algebra of the form $(B,C,B)$.
\end{rem}

\begin{nota}\label{compregn}
\rm Let $(X,\tau)$ be a topological space. We denote by
$CR(X,\tau)$ the family of all compact regular closed subsets of
$(X,\tau)$. We will often write  $CR(X)$ instead of $CR(X,\tau)$.
\end{nota}

\begin{fact}\label{stanlocn}{\rm (\cite{R})}
Let $(X,\tau)$ be a locally compact Hausdorff space. Then the triple
$$(RC(X,\tau),\rho_{(X,\tau)}, CR(X,\tau))$$
  is a complete local contact algebra; it is called a
{\em standard local contact algebra}.
\end{fact}

\begin{defi}\label{C_rho}{\rm (\cite{39})}
\rm
Let $(B,\rho,\mathbb{B})$ be a local contact algebra. Define a binary relation $"C_{\rho}"$ on $B$ by

\begin{equation}\label{C_rho e}
aC_{\rho}b\iff (a\rho b\ \mbox{\rm or}\ a, b\not\in\mathbb{B}).
\end{equation}

\noindent It is called the Alexandroff extension of $\rho$ relatively to the LCA $(B,\rho,\mathbb{B})$ (or, when there is no ambiguity, simply, the Alexandroff extension of $\rho$).
\end{defi}

\begin{lm}\label{2.12}{\rm (\cite{39})}
Let $(B,\rho,\mathbb{B})$ be a local contact algebra. Then $(B,C_{\rho})$,
where $C_{\rho}$ is the Alexandroff extension of $\rho$, is a normal contact algebra.
\end{lm}

\begin{defi}\label{bounded}
\rm
Let $(B,\rho,\mathbb{B})$ be a local contact algebra. We will say
that $\sigma$ is a cluster in $(B,\rho,\mathbb{B})$ if $\sigma$ is a cluster in the NCA $(B,C_{\rho})$. A cluster $\sigma$ in $(B,\rho,\mathbb{B})$  is called bounded if $\sigma\cap\mathbb{B}\not=\emptyset$.
\end{defi}

\begin{lm}\label{2.14}{\rm \cite{39}}
Let $(B,\rho,\mathbb{B})$ be a local contact algebra and let $1\not\in\mathbb{B}$.
Then $\sigma_{\infty}^{(B,\rho,\mathbb{B})}=\{b\in B\st b\not\in\mathbb{B}\}$ is a cluster in $(B,\rho,\mathbb{B})$. (Sometimes we will simply write $\sigma_{\infty}$ instead of $\sigma_{\infty}^{(B,\rho,\mathbb{B})}$.)
\end{lm}

\begin{nota}
\rm
Let $(X,\tau)$ be a topological space. If $x\in X$ then we set:

\begin{equation}\label{sigma_x}
\sigma_x=\{F\in RC(X)\st x\in F\}.
\end{equation}

\noindent for every $x\in X$, $\sigma_x$ is a bounded cluster in the standard local contact algebra $(RC(X,\tau),\rho_{(X,\tau)}, CR(X,\tau))$.
\end{nota}

The next  theorem was proved by Roeper \cite{R} (but its
particular case concerning compact Hausdorff spaces and NC-algebras was
proved by de Vries \cite{deV}).

\begin{theorem}\label{roeperl}{\rm (P. Roeper \cite{R}
for locally compact spaces and de Vries \cite{deV} for compact
spaces)}
There exists a bijective correspondence $\Psi^t$ between the
class of all (up
to homeomorphism) locally compact Hausdorff spaces and the class of all (up to isomorphism) CLC-algebras; its
restriction to the class of all (up to homeomorphism) compact Hausdorff spaces  gives a
bijective correspondence between the later class and the class of
all (up to isomorphism) CNC-algebras.
\end{theorem}

We will now recall (following \cite{39}) the definition of the correspondence $\Psi^t$ (mentioned
in the above theorem) and some other facts and notation which
 will be used later on.

Let $(X,\tau)$ be a locally compact Hausdorff space. Set

\begin{equation}\label{Psi^t}
\Psi^t(X,\tau) = (RC(X,\tau),\rho_{(X,\tau)},CR(X,\tau))
\end{equation}

Let $\unl{B}_{\, l}=(B,\rho,\BBBB)$ be a complete
local contact algebra. Let $C=C_\rho$ be the Alexandroff extension
of $\rho$. Then $(B,C)$ is
a complete normal contact algebra. Put $X=\CL(B,C)$ and let $\TT$
be the topology on $X$ having as a closed base the family
$\{\l_{(B,C)}(a)\st a\in B\}$ where, for every $a\in B$,
$\l_{(B,C)}(a) = \{\s \in X\st  a \in \s\}$.
Sometimes we will write simply $\l_B$ instead of $\l_{(B,C)}$.
 Note that
$X\stm \l_B(a)= \int(\l_B(a^*))$,
the family
$\{\int(\l_B(a))\st a\in B\}$
is an open
base of
$(X,\TT)$
and, for every $a\in B$,
$\l_B(a)\in RC(X,\TT)$.
It can be proved that
$\l_B:(B,C)\lra (RC(X),\rho_X)$
is a CA-isomorphism.
Further,
$(X,\TT)$
is a compact Hausdorff space.

 Let $1\in\BBBB$. Then $C=\rho$ and $\BBBB=B$, so
that $(B,\rho,\BBBB)=(B,C,B)=(B,C)$ is a complete normal contact
algebra, and we put
\begin{equation}\label{phiapcn}
\Psi^a(B,\rho,\BBBB)=\Psi^a(B,C,B)=\Psi^a(B,C)=(X,\TT).
\end{equation}

\medskip

 Let $1\not\in\BBBB$. Then
the set $\s_\infty=\{b\in B\st b\not\in\BBBB\}$ is a cluster in
$(B,C)$ and, hence, $\s_\infty\in X$.  Let $L=X\stm\{\s_\infty\}$.
Then
$L=\BClu(B,\rho,\BBBB)$,
i.e.
$L$
is the set of all
bounded clusters of
$ (B,C_\rho)$
(sometimes we will write $L_{\unl{B}_{\, l}}$ or $L_B$ instead of
$L$);
 let the topology $\tau(=\tau_{\unl{B}_{\, l}})$ on $L$ be the
subspace topology, i.e. $\tau=\TT|_L $. Then $(L,\tau)$ is a
locally compact Hausdorff space. We put
\begin{equation}\label {phiapc}
\Psi^a(B,\rho,\BBBB)=(L,\tau).
\end{equation}

Let
$\l^l_{\unl{B}_{\, l}}(a)=\l_{(B,C_\rho)}(a)\cap L$,
for each $a\in B$. We will write simply $\l^l_B$ (or even
$\l_{(A,\rho,\BBBB)}$ when $\BBBB\neq A$) instead of
$\l^l_{\unl{B}_{\, l}}$ when this does not lead to ambiguity. One
can show that:

\noindent (I) $L$ is a dense subset of $X$;\\
(II) $\l^l_B$ is a Boolean isomorphism of the Boolean algebra $B$ onto
the
Boolean algebra $RC(L,\tau)$;\\
(III) $b\in\BBBB$ iff $\l^l_B(b)\in CR(L)$;\\
(IV) $a\rho b$ iff $\l^l_B(a)\cap \l^l_B(b)\neq\ems$.\\
\noindent Hence, $X=\a L$  and
$\l^l_B: (B,\rho,\BBBB)\lra (RC(L),\rho_L, CR(L))$
is an
LCA-isomorphism.
 Note also that for every $b\in B$,
$\int_{L_B}(\l^l_B(b))=L_B\cap\int_X(\l_B(b))$.

 For every CLCA $(B,\rho,\BBBB)$ and every $a\in B$,
set
\begin{equation}\label{lbg}
\l^g_{\unl{B}_{\,l}}(a)=\l_{(B,C_\rho)}(a)\cap\Psi^a(B,\rho,\BBBB).
\end{equation}
We will write simply $\l^g_B$ instead of $\l^g_{\unl{B}_{\, l}}$
when this does not lead to ambiguity.
 Thus, when $1\in\BBBB$,
we have that $\l^g_B=\l_B$, and  if $1\nin\BBBB$  then
$\l^g_B=\l^l_B$. Hence we
get that
\begin{equation}\label{hapisomn}
\l^g_B: (B,\rho,\BBBB)\lra (\Psi^t\circ\Psi^a)(B,\rho,\BBBB)
\mbox{ is an LCA-isomorphism.}
\end{equation}

We have that:
\begin{equation}\label{eel}
\mbox{the family } \{\int_{\Psi^a(B,\rho,\BBBB)}(\l_B^g(a))\st
a\in \BBBB\} \mbox{ is an open base of } \Psi^a(B,\rho,\BBBB).
\end{equation}

Let  $(L,\tau)$ be a  locally compact Hausdorff space,
 $B=RC(L,\tau)$, $\BBBB=CR(L,\tau)$ and $\rho=\rho_L$. Then
$(B,\rho,\BBBB)=\Psi^t(L,\tau)$. It
can be shown that the map
\begin{equation}\label{homeo}
t_{(L,\tau)}:(L,\tau)\lra\Psi^a(\Psi^t(L,\tau)),
\end{equation}
defined by  $t_{(L,\tau)}(x)=\{F\in RC(L,\tau)\st x\in
F\}(=\s_x)$, for all $x\in L$, is a homeomorphism.

Therefore $\Psi^a(\Psi^t(L,\tau))$ is homeomorphic to $(L,\tau)$
and $\Psi^t(\Psi^a(B,\rho,\BBBB))$ is LCA-isomorphic to
$(B,\rho,\BBBB)$.

Note that  if $(A,\rho,\mathbb{B})$ is an LCA, $X=\Psi^a(A,\rho,\mathbb{B})$ and $(B,\eta,\mathbb{B}')=\lambda^g_B(A,\rho,\mathbb{B})$ then for every $a\in RC(X)$, $a=\bigvee\{b\in\mathbb{B}'\st b\ll_{\rho_X}a\}$ holds. Hence, for every $a\in A$,

\begin{equation}\label{21}
a=\bigvee\{b\in\mathbb{B}\st b\ll_{\rho}a\}.
\end{equation}
%%%%%%%%%%%%%%%%%%%%%%%%%%%%%%%%%%%%%%%%%%%%%%%%%%%%%%%%%%%%%%%%%%%%%%%%%%%%%

\begin{defi}\label{lideal}{\rm (\cite{D-AMH1-10})}
\rm Let $(A,\rho,\mathbb{B})$ be an LCA. An ideal $I$ of $A$ is called
a {\em $\delta$-ideal} if $I\subseteq \mathbb{B}$ and for any $a\in I$ there
exists $b\in I$ such that $a\ll_\rho b$. If $I_1$ and $I_2$ are
two $\delta$-ideals of $(A,\rho,\mathbb{B})$ then we put $I_1\le I_2$ iff
$I_1\subseteq I_2$. We will denote by $(I(A,\rho,\mathbb{B}),\le)$ the poset
of all $\delta$-ideals of $(A,\rho,\mathbb{B})$.
\end{defi}

\begin{fact}\label{dideal}{\rm (\cite{D-AMH1-10})}
Let $(A,\rho,\mathbb{B})$ be an LCA. Then, for every $a\in A$, the set
$I_a=\{b\in\mathbb{B}\st b\ll_\rho a\}$ is a $\delta$-ideal. Such
$\delta$-ideals will be called\/ {\em principal $\delta$-ideals}.
\end{fact}

Recall that a {\em frame} is a complete lattice $L$ satisfying the
infinite distributive law $a\wedge\bigvee S=\bigvee\{a\wedge s\st s\in
S\}$, for every $a\in L$ and every $S\subseteq L$.

\begin{fact}\label{frlid}{\rm (\cite{D-AMH1-10})}
Let $(A,\rho,\mathbb{B})$ be an LCA. Then the poset
$(I(A,\rho,\mathbb{B}),\le)$ of all $\delta$-ideals of $(A,\rho,\mathbb{B})$
is a frame. The finite meets and arbitrary joins in $I(A,\rho,\BBBB)$ coincide with the corresponding
operations in the frame $Idl(A)$ of all ideals of $A$.
\end{fact}

We will often use the following elementary fact: the join $\bv\{I_\g\st\g\in\GA\}$ of a family of ideals of a distributive lattice $A$ in the frame $Idl(A)$ of all ideals of $A$ is
the  set $I=\{\bv\{x_\g\st\g\in\GA_1\}\st \GA_1\sbe\GA,\GA_1$ is finite, $x_\g\in I_\g$ for every $\g\in\GA_1\}$ of elements of $A$ (see, e.g., \cite{Dw}).

\begin{pro}\label{clustuniq}{\rm (\cite{D-AMH1-10})}
Let $\sigma_1,\sigma_2\in\Psi^a(A,\rho,\mathbb{B})$, where $(A,\rho,\BBBB)$ is a CLCA, and
$\sigma_1\cap\mathbb{B}=\sigma_2\cap\mathbb{B}$. Then $\sigma_1=\sigma_2$.
\end{pro}

Recall that if $A$ is a distributive lattice then an element $p\in A\stm\{1\}$
is called a {\em prime element of}\/ $A$ if for each $a,b\in A$,
$a\we b= p$ implies that $a= p$ or $b=p$. The prime elements of the frame $I(A,\rho,\BBBB)$,
where $(A,\rho,\BBBB)$ is an LCA, will be called {\em prime $\d$-ideals of} $(A,\rho,\BBBB)$.

\begin{pro}\label{clustideal}{\rm (\cite{D-AMH1-10})}
Let $(A,\rho,\mathbb{B})$ be a CLCA. If  $\sigma\in\Psi^a(A,\rho,\mathbb{B})$
then $\mathbb{B}\setminus\sigma=J_{\sigma}$ is a prime $\d$-ideal of $(A,\rho,\BBBB)$. If $J$ is a prime
$\d$-ideal of $(A,\rho,\BBBB)$ then there exists a unique $\s\in\Psi^a(A,\rho,\BBBB)$ such that $\s\cap\BBBB=\BBBB\stm J$.
\end{pro}

\begin{theorem}\label{iota}{\rm (\cite{D-AMH1-10})}
Let $(A,\rho,\mathbb{B})$ be a CLCA, $X=\Psi^a(A,\rho,\mathbb{B})$ and ${\cal O}(X)$ be the frame of all open subsets of $X$. Then there exists a frame isomorphism
$$\iota:(I(A,\rho,\mathbb{B}),\leq)\longrightarrow({\cal O}(X),\subseteq),\ \  I\mapsto\bigcup\{\lag(a)\st a\in I\},$$
where $(I(A,\rho,\mathbb{B}),\leq)$ is the frame of all $\delta$-ideals of $(A,\rho,\mathbb{B})$.
\end{theorem}

%----------------------------------------------------------------------------------
\section{A new duality theorem}
%----------------------------------------------------------------------------------

\begin{nota}
\rm
We denote by {\bf HLC} the category of all locally compact Hausdorff spaces and all continuous mappings between them.
\end{nota}

\begin{defi}\label{maincat}
\rm
Let {\bf MDHLC} be the category whose objects are all CLCAs and whose morphisms $\p : (A,\rho,\mathbb{B})\longrightarrow(B,\eta,\mathbb{B}')$ are all multi-valued maps which satisfy the following conditions:
\smallskip

\noindent (M1) For every $a\in A$, $\p (a)\in I(B,\eta,\mathbb{B}')$;

\noindent (M2) $\p (a\wedge b)=\p (a)\we \p (b)$, for every $a,b\in A$;

\noindent (M3) $\p (a)= \bv\{\p (b)\st b\in\mathbb{B}, b\ll a\}$, for every $a\in A$;

\noindent (M4) $\p (0)=\{0\}$;

\noindent (M5) If  $a_i,b_i\in\mathbb{B}$, $a_i\ll b_i$, where $i=1,2$, then $\p (a_1\vee a_2)\subseteq \p (b_1)\vee \p (b_2)$;

\noindent (M6) For every $b\in\mathbb{B}'$, there exists an $a\in\mathbb{B}$ such that $b\in \p (a)$.
\smallskip

The composition $\diamond$ between two morphisms $\p :(A_1,\rho_1,\mathbb{B}_1)\longrightarrow(A_2,\rho_2,\mathbb{B}_2)$ and $\psi:(A_2,\rho_2,\mathbb{B}_2)\longrightarrow (A_3,\rho_3,\mathbb{B}_3)$ is defined by $(\psi\diamond \p )(a)=\bigvee\{\psi(b)\st b\in \p (a)\}$. The identity morphism $i_A:(A,\rho,\mathbb{B})\longrightarrow (A,\rho,\mathbb{B})$ is defined by $i_A(a)=I_a$ (see \ref{dideal} for $I_a$).
\end{defi}

\begin{rem}\label{remmdhlc}
\rm
Using Fact \ref{frlid}, it can be easily seen that in the axiom (M2) the expression $``\p (a)\we \p (b)$" can be replaced by $``\p (a)\cap \p (b)$", and,
in (M3), $``\bv$" can be replaced by $``\bigcup$". Note also that the expression $``\bigvee\{\psi(b)\st b\in \p (a)\}$"   can be written down
in the form $``\bv\psi(\p(a))$", and hence
$(\psi\diamond \p )(a)=\bigvee\psi(\p(a))$, i.e.  our definition of the composition between two morphisms in the category $\MDHLC$ is enough natural.
\end{rem}

\begin{pro}\label{mdhlccat}
{\bf MDHLC} is a category.
\end{pro}

\doc
We will first  prove that for every $(A,\rho,\mathbb{B})$, $i_A$ is an {\bf MDHLC}-morphism.
Indeed, it is obvious that (M1), (M2) and (M4) are satisfied. Since (BC1) implies that
$i_A(a)=I_a=\bigvee\{I_b\st b\in I_a\}$, we get that (M3) is fulfilled. We will now show that condition (M5) is fulfilled.
Let $a_i,b_i\in\mathbb{B}$, $a_i\ll b_i$, $i=1,2$. We have to show that
$I_{a_1\vee a_2}\subseteq I_{b_1}\vee I_{b_2}$. Let $c\ll a_1\vee a_2$. Then $c=(c\wedge a_1)\vee(c\wedge a_2)$.
Since $c\wedge a_1\le a_1\ll b_1$ and $c\wedge a_2\le a_2\ll b_2$, we get that
$c\wedge a_1\in I_{b_1}$ and $c\wedge a_2\in I_{b_2}$.
Hence $c=(c\wedge a_1)\vee(c\wedge a_2)\in I_{b_1}\vee I_{b_2}$. So, $I_{a_1\vee a_2}\subseteq I_{b_1}\vee I_{b_2}$.
For verifying (M6), let $b\in\mathbb{B}$; then, by (BC1), there exists an $a\in\mathbb{B}$ such that $b\ll a$;
hence $b\in I_a=i_A(a)$. So, $i_A$ is a {\bf MDHLC}-morphism.

Let $\p _1:(A_1,\rho_1,\mathbb{B}_1)\longrightarrow (A_2,\rho_2,\mathbb{B}_2)$ and
$\p _2:(A_2,\rho_2,\mathbb{B}_2)\longrightarrow (A_3,\rho_3,\mathbb{B}_3)$ be {\bf MDHLC}-morphisms.
We will prove that $\p =\p _2\diamond \p _1$ is an {\bf MDHLC}-morphism.
We have that $\p (a)=\bigvee\{\p _2(b)\st b\in \p _1(a)\}$. The axiom (M1) is obviously fulfilled.
Further, for every $a_1, a_2\in A_1$,
$$\p (a_1\wedge a_2)=\bigvee\{\p _2(b)\st b\in \p _1(a_1\wedge a_2)\}= \bigvee\{\p _2(b)\st b\in \p _1(a_1)\cap \p _1(a_2)\}$$
 and

$\begin{array}{rl}
\p (a_1)\wedge \p (a_2)&=\bigvee\{\p _2(b_1)\st b_1\in \p _1(a_1)\}\wedge\bigvee\{\p _2(b_2)\st b_2\in \p _1(a_2)\}\\
&=\bigvee\{\p _2(b_2)\wedge\bigvee\{\p _2(b_1)\st b_1\in \p _1(a_1)\}\st b_2\in \p _1(a_2)\}\\
&=\bigvee\{\bigvee\{\p _2(b_1)\wedge \p _2(b_2)\st b_1\in \p _1(a_1)\}\st b_2\in \p _1(a_2)\}\\
&=\bigvee\{\p _2(b_1\wedge b_2)\st b_1\in \p _1(a_1), b_2\in \p _1(a_2)\}.
\end{array}$

If, for $i=1,2$, $b_i\in \p _1(a_i)$ then $b_1\wedge b_2=b\in \p _1(a_1)\cap \p _1(a_2)$.
So, $\p (a_1)\wedge \p (a_2)\subseteq \p (a_1\wedge a_2)$.
Conversely, from $b\in \p _1(a_1)\cap \p _1(a_2)$ and $b=b\wedge b$, we get that
 $\p (a_1\wedge a_2)\subseteq \p (a_1)\wedge \p (a_2)$. Hence, condition (M2) is satisfied.

We will prove that $\p (a)=\bv\{\p(b)\st b\in\BBBB, b\ll a\}$ for every $a\in A$, i.e.
$\bv\{\p_2(c)\st c\in\p_1(a)\}=\bigvee\{\p _2(d)\st d\in \p _1(b), b\in\mathbb{B}, b\ll a\}$.
Let $c\in \p _1(a)$. Then, by (M3) and Remark \ref{remmdhlc}, there exists  $b\in\mathbb{B}$  such that $b\ll a$ and $c\in \p _1(b)$.
Conversely, let $d\in \p _1(b)$, $b\in\mathbb{B}$, $b\ll a$. Then $d\in \p _1(a)$. Hence, the axiom (M3) is fulfilled.

Since $\p (0)=\bigvee\{\p _2(b)\st b\in \p _1(0)\}=\bigvee\{\p _2(b)\st b\in\{0\}\}=\p _2(0)=\{0\}$, we get that condition (M4) is satisfied.

Let $a_i, b_i\in\mathbb{B}$, $a_i\ll b_i$, $i=1, 2$. We will prove that
$\p (a_1\vee a_2)\subseteq \p (b_1)\vee \p (b_2)$, i.e.
$$\bigvee\{\p _2(c)\st c\in \p _1(a_1\vee a_2)\}\subseteq\bigvee\{\p _2(d)\st d\in \p _1(b_1)\}\vee\bigvee\{\p _2(e)\st e\in \p _1(b_2)\}.$$
Let $c\in \p _1(a_1\vee a_2)$. Then $c\in \p _1(b_1)\vee \p _1(b_2)$, i.e.
there exist $d_1\in \p _1(b_1)$ and $e_1\in \p _1(b_2)$ such that $c=d_1\vee e_1$.
There exists  $d\in \p _1(b_1)$ such that $d_1\ll d$ and there exists  $e\in \p _1(b_2)$ such that $e_1\ll e$.
Then $\p _2(c)=\p _2(d_1\vee e_1)\subseteq \p _2(d)\vee \p _2(e)$. So, the axiom (M5) is satisfied.

Let $c\in\mathbb{B}_3$. Then there exists  $b\in\mathbb{B}_2$ such that $c\in \p _2(b)$.
There exists  $a\in\mathbb{B}_1$ such that $b\in \p _1(a)$. Hence $c\in \p (a)$. So, condition (M6) is also fulfilled.

Hence $\p _2\diamond \p _1$ is an {\bf MDHLC}-morphism.

We will now show that the composition
is associative.
Let $\p :(A_1,\rho_1,\mathbb{B}_1)\longrightarrow(A_2,\rho_2,\mathbb{B}_2)$, $\psi :(A_2,\rho_2,\mathbb{B}_2)\longrightarrow(A_3,\rho_3,\mathbb{B}_3)$ and $\chi:(A_3,\rho_3,\mathbb{B}_3)\longrightarrow(A_4,\rho_4,\mathbb{B}_4)$ be {\bf MDHLC}-morphisms. We have that, for every $a\in A_3$,

$\begin{array}{rl}
(\p \diamond(\psi \diamond \chi))(a)&=\bigvee\{\p (b)\st b\in(\psi \diamond \chi)(a)\}\\
&=\bigvee\{\p (b)\st b\in\bigvee\{\psi (c)\st c\in \chi(a)\}\},
\end{array}$

\noindent and

$\begin{array}{rl}
((\p \diamond \psi )\diamond \chi)(a)&=\bigvee\{(\p \diamond \psi )(c)\st c\in \chi(a)\}\\
&=\bigvee\{\bigvee\{\p (b)\st b\in \psi (c)\}\st c\in \chi(a)\}\\
&=\bv\{\p(b)\st b\in\bigcup\{\psi(c)\st c\in\chi(a)\}\}.
\end{array}$

\noindent Let $b\in\bv\{\psi(c)\st c\in\chi(a)\}$. Then $b=\bv\{b_i\st i\in\{1,\ldots,n\}\}$, for some $n\in\NNNN^+$, where, for every
$i\in\{1,\ldots,n\}$, $b_i\in\psi(c_i)$ and $c_i\in\chi(a)$. Setting $c=\bv\{c_i\st i\in\{1,\ldots,n\}\}$, we get that $c\in\chi(a)$ and,
by (M1) and (M2), $\bv\{\psi(c_i)\st i\in\{1,\ldots,n\}\}\sbe\psi(c)$. Therefore $b\in\psi(c)$. We get that
$\bv\{\psi(c)\st c\in\chi(a)\}=\bigcup\{\psi(c)\st c\in\chi(a)\}$. Hence, the composition $``\diamond$" is associative.

Finally, if $\p:(A,\rho,\BBBB)\lra (B,\eta,\BBBB\ap)$ is a {\bf MDHLC}-morphism then, for every $a\in A$,
$(\p\diamond i_A)(a)=\bv\{\p(b)\st b\in I_a\}=\bv\{\p(b)\st b\in\BBBB, b\ll a\}=\p(a)$ (since $\p$ satisfies condition (M3)), and
$(i_B\diamond\p)(a)=\bv\{I_b\st b\in\p(a)\}=\p(a)$. Hence, $\p\diamond i_A=\p$ and $i_B\diamond\p=\p$.

All this shows that {\bf MDHLC} is a category.
\sqs

\begin{pro}\label{Delta^t}
Let $f:X\longrightarrow Y$ be an {\bf HLC}-morphism. Define a map $\varphi_f:\Psi^t(Y)
\longrightarrow\Psi^t(X)$
by:
\begin{equation}\label{Delata^te}
\forall G\in RC(Y),\ \ \ \varphi_f(G)=\{F\in CR(X)\st F\subseteq f^{-1}(\int (G))\}.
\end{equation}
\noindent Then $\varphi_f$ is an {\bf MDHLC}-morphism.
\end{pro}

\doc We have to prove that $\varphi_f$ satisfies the conditions (M1)-(M6) from Definition \ref{maincat}. We start by showing that for each $G\in RC(Y)$, $\varphi_f(G)$ is a $\delta$-ideal. Obviously, $\varphi_f(G)$ is a lower set. If $F_1, F_2\in\varphi_f(G)$ then $F_1\vee F_2=F_1\cup F_2\in\varphi_f(G)$. So, $\varphi_f(G)$ is an ideal. If $F\in\varphi_f(G)$ then $F$ is compact and $F\subseteq f^{-1}(\int (G))$. Hence there exists an open $U\subseteq X$ such that $\cl(U)$ is compact and $F\subseteq U\subseteq \cl(U)\subseteq f^{-1}(\int (G))$. Then $\cl(U)\in CR(X)$ and hence $\cl(U)\in\varphi_f(G)$. So, $\varphi_f(G)$ is a $\delta$-ideal. Thus, condition (M1) is fulfilled.

Let $G,H\in RC(Y)$. Then $$\varphi_f(G\wedge H)=\{F\in CR(X)\st F\subseteq f^{-1}(\int(G\wedge H))\}$$ and
$$\begin{array}{rl}
\varphi_f(G)\cap\varphi_f(H)&=\{F\in CR(X)\st F\subseteq f^{-1}(\int (G)),\ F\subseteq f^{-1}(\int (H))\}\\
&=\{F\in CR(X)\st F\subseteq f^{-1}(\int(G\cap H))\}.
\end{array}$$
\noindent Since $\int(G\cap H)$ is a regular open set, we get that $\int(G\wedge H)=\int(\cl(\int(G\cap H)))=\int(G\cap H)$. So, $\varphi_f(G\wedge H)=\varphi_f(G)\cap\varphi_f(H)$. Thus, the axiom (M2) is satisfied.

For verifying (M3), we have to prove that
$\{F\in CR(X)\st F\subseteq f^{-1}(\int (G))\} =$ $\bigvee\{\{F'\in CR(X)\st F'\subseteq f^{-1}(\int (H))\}\st H\in CR(Y),\ H\subseteq \int (G)\}$.
It is obvious that the right part is a subset of the left part. For proving the converse inclusion,
let $F\in CR(X)$ and $F\subseteq f^{-1}(\int (G))$. Then $f(F)\subseteq \int (G)$ and $f(F)$ is compact. Let $\Omega=\{\int (H)\st H\in CR(Y),\ H\subseteq \int (G)\}$. Then $\bigcup\Omega=\int (G)$. Hence $\Omega$ covers $f(F)$. Therefore there exist $H_1,\ldots ,H_n$ such that $\int (H_1),\ldots,\int (H_n)\in\Omega$ and $\displaystyle f(F)\subseteq\bigcup_{i=1}^n \int (H_i)\subseteq\bigcup_{i=1}^n H_i\subseteq \int (G)$. Set $\displaystyle H=\bigcup_{i=1}^n H_i$. Then $H\in CR(Y)$ and $H\subseteq \int (G)$.
Since $\displaystyle \bigcup_{i=1}^n \int (H_i)\subseteq \int(\bigcup_{i=1}^n H_i)$, we get that $f(F)\subseteq \int (H)$, i.e. $F\sbe f\inv(\int (H))$.
Thus, condition (M3) is fulfilled.

We have that $0=\emptyset$, so $\varphi_f(\emptyset)=\{F\in CR(X)\st F\subseteq f^{-1}\{\emptyset\}\}=\{\emptyset\}=I_{\emptyset}$. Therefore, $\varphi_f$ satisfies  condition (M4).

For verifying the axiom (M5), we have to prove that for every $G_i, H_i\in CR(Y)$ such that $G_i\subseteq \int (H_i)$, where $i=1,2$, the following inclusion holds:
 $$\{F\in CR(X)\st F\subseteq f^{-1}(\int(G_1\cup G_2))\}\subseteq$$ $$\{F'\in CR(X)\st F'\subseteq f^{-1}(\int (H_1))\}\vee\{F''\in CR(X)\st F''\subseteq f^{-1}(\int (H_2))\}.$$ Let $F\in CR(X)$ and $F\subseteq f^{-1}(\int(G_1\cup G_2))$. Then $$F\subseteq f^{-1}(G_1\cup G_2)=f^{-1}(G_1)\cup f^{-1}(G_2)\subseteq f^{-1}(\int (H_1))\cup f^{-1}(\int (H_2)).$$ Obviously, $\Omega_i=\{\int (K)\st K\in CR(X),\ K\subseteq f^{-1}(\int (H_i))\}$ covers $f^{-1}(\int (H_i))$, for $i=1,2$. Then $\Omega=\Omega_1\cup\Omega_2$ is a cover of $f^{-1}(\int (H_1))\cup f^{-1}(\int (H_2))$ and hence $F\subseteq\bigcup\Omega$. Since $F$ is compact, there exist
$\int (K_1),\ldots,\int (K_m)\in\Omega _1$ and $\int (K_1'),\ldots,\int (K_n')\in\Omega_2$ such that
$\displaystyle F\subseteq\bigcup_{i=1}^m \int (K_i)\cup\bigcup_{j=1}^n \int (K_j')$. Put
$\displaystyle F_1=\bigcup_{i=1}^m
K_i$ and $\displaystyle F_2=\bigcup_{j=1}^n
K_j'$. Then $F_i\in CR(X)$ and $F_i\subseteq f^{-1}(\int (H_i))$,  where $i=1,2$. Therefore $F\subseteq F_1\cup F_2$ and $F_1\cup F_2\in\varphi_f(H_1)\vee\varphi_f(H_2)$. Hence $F\in\varphi_f(H_1)\vee\varphi_f(H_2)$.

Finally, we will show that (M6) is fulfilled. Let $F\in CR(X)$. For every $y\in f(F)$ there exists a neighborhood $O_y$ of $y$ such that $\cl(O_y)$ is compact.
 Since $f(F)$ is  compact,  there exist $y_1,\ldots, y_n\in f(F)$ such that $\displaystyle f(F)\subseteq\bigcup_{i=1}^n O_{y_i}$. Let $\displaystyle G=\bigcup_{i=1}^n \cl(O_{y_i})$. Then $G\in CR(Y)$ and $f(F)\subseteq \int (G)$. Hence $F\subseteq f^{-1}(\int (G))$, i.e. $F\in\varphi_f(G)$. \sqs

\begin{pro}\label{functa}
For each $X\in|{\bf HLC}|$, set $\Delta^t(X)=\Psi^t(X)$  (see Theorem \ref{roeperl} for the notation $\Psi^t$), and for each $f\in{\bf HLC}(X,Y)$, put $\Delta^t(f)=\varphi_f$ (see Proposition \ref{Delta^t} for the notation $\p_f$). Then $\DE^t:\HLC\lra\MDHLC$ is a contravariant functor.
\end{pro}

\doc Let $X\in|\HLC|$ and $(A,\rho,\BBBB)=\DE^t(X)$. We will show that $\DE^t(id_X)=i_A$. Indeed, let $\p=\DE^t(id_X)$. Then, by (\ref{Delata^te}),
$\p(G)=\{F\in CR(X)\st F\sbe\int(G)\}=\{a\in\BBBB\st a\ll G\}=I_G=i_A(G)$, for every $G\in RC(X)$ $(=A)$. Thus $\DE^t(id_X)=i_A$.

Let now $f_1\in\HLC(X_1,X_2)$,  $f_2\in\HLC(X_2,X_3)$ and $f=f_2\circ f_1$. We will show that $\DE^t(f)=\DE^t(f_1)\diamond\DE^t(f_2)$. Set, for short,
$\p=\DE^t(f)$,  $\p_1=\DE^t(f_1)$ and $\p_2=\DE^t(f_2)$. Then, for every $G_3\in RC(X_3)$, we have that
$\p_2(G_3)=\{F_2\in CR(X_2)\st f_2(F_2)\sbe \int(G_3)\},$
\begin{equation}\label{eq1}
\p(G_3)=\{F_1\in CR(X_1)\st F_1\sbe f_1\inv(f_2\inv(\int(G_3)))\}
\end{equation}
  \noindent and
  \medskip

\noindent$\begin{array}{rl}
(\p_1\di\p_2)(G_3)&= \bv\{\p_1(F_2)\st F_2\in\p_2(G_3)\}\\
&= \bv\{\{F_1\in CR(X_1)\st f_1(F_1)\sbe \int(F_2)\}\st\\
&\ \ \ \ \ \ \ \ \ F_2\in CR(X_2), f_2(F_2)\sbe \int(G_3)\}\\
&= \{\bigcup\{ F_1^i\st i=1,\ldots,k\}\st k\in\NNNN^+,(\fa i=1,\ldots,k)[(F_1^i\in CR(X_1))\we\\
&\ \ \ \ \ \ \ \ \ ((\ex F_2^i\in CR(X_2))(f_1(F_1^i)\sbe\int(F_2^i)\sbe F_2^i\sbe f_2\inv(\int(G_3))))]\}\\
&=\{F_1\in CR(X_1)\st (\ex F_2\in CR(X_2))\\
&\ \ \ \ \ \ \ \ \ (f_1(F_1)\sbe\int(F_2)\sbe F_2\sbe f_2\inv(\int(G_3)))\}.
\end{array}$

\smallskip

\noindent We have to show that $\p(G_3)=(\p_1\di\p_2)(G_3)$, i.e. that the corresponding right sides $R$ and $R_{1,2}$ of (\ref{eq1}) and the equation after it are equal.
Let $F_1\in R$. Then $F_1\in CR(X_1)$ and $f_1(F_1)\sbe f_2\inv(\int(G_3))$. Since $f_1(F_1)$ is a compact subset of $X_2$, there exists $F_2\in CR(X_2)$ such that
$f_1(F_1)\sbe\int(F_2)\sbe F_2\sbe f_2\inv(\int(G_3))$. Thus, $F_1\in R_{1,2}$. Conversely, if $F_1\in R_{1,2}$ then $F_1\in CR(X_1)$ and
there exists $F_2\in CR(X_2)$ such that
$f_1(F_1)\sbe\int(F_2)\sbe F_2\sbe f_2\inv(\int(G_3))$. Then $F_1\sbe f_1\inv(F_2)\sbe f_1\inv(f_2\inv(\int(G_3)))$. Therefore $F_1\in R$.
So, we have proved that $\DE^t(f)=\DE^t(f_1)\diamond\DE^t(f_2)$. All this shows that $\DE^t$ is a contravariant functor. \sqs

\begin{pro}\label{Delta^a}
Let $\varphi:(A,\rho,\mathbb{B})\longrightarrow(B,\rho',\mathbb{B}')$ be an {\bf MDHLC}-morphism. Define a map $f_{\varphi}:\Psi^a(B,\rho',\mathbb{B}')\longrightarrow\Psi^a(A,\rho,\mathbb{B})$ by setting

\begin{equation}\label{Delata^ae}
\forall \sigma'\in\Psi^a(B,\rho',\mathbb{B}'), f_{\varphi}(\sigma')\cap\mathbb{B}=\{a\in\mathbb{B}\st (\forall b\in A)((a\ll_\rho b)\rightarrow (\varphi(b)\cap\sigma'\not=\emptyset))\}.
\end{equation}
\noindent Then $f_{\varphi}$ is defined correctly and $f_\p$ is an {\bf HLC}-morphism.
\end{pro}

\doc Let $\sigma'\in\Psi^a(B,\rho',\mathbb{B}')$. Set $J=\BBBB\stm (f_{\varphi}(\sigma')\cap\BBBB)$. We will  first prove that $J$  is a prime $\delta$-ideal
of $(A,\rho,\BBBB)$. Note that $J=\{a\in\mathbb{B}\st\exists b\in\mathbb{B}$ such that $a\ll_{\rho}b$ and $\varphi(b)\cap\sigma'=\emptyset\}$.

Obviously, $J $ is a lower set. By (M4), $0\in J $ (because $0\ll_{\rho}0$). Let $a,b\in J $. Then there exist $a',b'\in\mathbb{B}$ such that $a\ll_{\rho}a'$, $b\ll_{\rho}b'$ and $\varphi(a')\cap\sigma'=\emptyset$, $\varphi(b')\cap\sigma'=\emptyset$. There exist $a'',b''\in\mathbb{B}$ such that $a\ll_{\rho}a''\ll_{\rho}a'$ and $b\ll_{\rho}b''\ll_{\rho}b'$. Hence, by (M5), $\varphi(a''\vee b'')\subseteq\varphi(a')\vee\varphi(b')$. Since, by Proposition \ref{clustideal},
 $\mathbb{B}\ap\setminus\sigma'$ is a $\delta$-ideal and $\varphi(a')\cup\varphi(b')\subseteq\mathbb{B}\ap\setminus\sigma'$, we get that $\varphi(a')\vee\varphi(b')\subseteq\mathbb{B}\ap\setminus\sigma'$. Thus $\varphi(a''\vee b'')\cap\sigma'=\emptyset$. Since $a\vee b\ll_{\rho}a''\vee b''$,
 we obtain that $a\vee b\in J $. Hence $J $ is an ideal.

Let $a\in J $. Then there exists  $b\in\mathbb{B}$ such that $a\ll_{\rho}b$ and $\varphi(b)\cap\sigma'=\emptyset$. There exists  $c\in\mathbb{B}$ such that $a\ll_{\rho}c\ll_{\rho}b$. Then, obviously, $c\in J $ and $a\ll_{\rho}c$. Hence $J $ is a $\delta$-ideal.

Let $I_1,I_2\in I(A,\rho,\BBBB)$ and $I_1\cap I_2=J $. Suppose that $J \neq I_i$, for $i=1,2$.
Hence there exists $a_i\in I_i\setminus J $, for $i=1,2$.
Then, for every $b\in\mathbb{B}$ such that $a_1\ll b$ or $a_2\ll b$, we have that $\varphi(b)\cap\sigma'\not=\emptyset$.
There exists $b_i\in I_i$ such that $a_i\ll b_i$, for $i=1,2$. Then $b_i\not\in J $, for $i=1,2$. Let $b=b_1\wedge b_2$. Then $b\in I_1\cap I_2=J $ and thus $\p(b)\cap\s\ap=\ems$. Using (M2), we get that $\varphi(b_1)\cap\varphi(b_2)\cap\sigma'=\emptyset$. There exists $d_i\in\varphi(b_i)\cap\sigma'$, for $i=1,2$. Since $\p(b_i)$ is a $\d$-ideal, there exists $l_i\in\varphi(b_i)$ such that $d_i\ll l_i$, for $i=1,2$. Then $l_i\in\sigma'$  but $l_i^*\not\in\sigma'$ (since $d_i(-C_{\rho})l_i^*$), where $i=1,2$. Hence $l_1^*\vee l_2^*\not\in\sigma'$. Then $l_1\wedge l_2\in\sigma'$. Moreover, $l_1\wedge l_2\in\varphi(b_1)\cap\varphi(b_2)\cap\sigma'$, which is a contradiction.

So, $J $ is a prime $\delta$-ideal. Obviously, $\BBBB\stm J =f_\p(\s')\cap\BBBB$. Now, by Proposition \ref{clustideal}, there exists a unique bounded cluster $\s$ in $(A,\rho,\BBBB)$ whose intersection with $\BBBB$ is equal to $\BBBB\stm J$. Thus $f_{\varphi}(\sigma')=\sigma$. All this shows that $f_\p$ is defined correctly.

We will now prove  that $f_{\varphi}$ is a continuous function. Let $F\in CR(X)$, where
$X=\Psi^a(A,\rho,\mathbb{B})$. Then there exists  $a\in\mathbb{B}$ such that  $F=\lambda_A^g(a)$.
Set $U= \int (F)$. Then $U=\int (\lambda_A^g(a))=X\setminus \lambda_A^g(a^*)$.
We will show that $f_{\varphi}^{-1}(U)=\iota_B(\varphi(a))$ ($=\bigcup\{\lambda_B^g(b)\st b\in\varphi(a)\}$).
Indeed, let $\sigma'\in f_{\varphi}^{-1}(U)$. Then $f_{\varphi}(\sigma')=\sigma\in U=X\setminus\lambda_A^g(a^*)$.
Hence $a^*\not\in\sigma$ and $a\in\sigma$. We have that

\begin{equation}\label{*}
\sigma\cap\mathbb{B}=\{c\in\mathbb{B}\st\forall d\in\mathbb{B}\mbox{ such that } c\ll d,\ \varphi(d)\cap\sigma'\not=\emptyset\}.
\end{equation}

\noindent We will prove   that $\varphi(a)\cap\sigma'\not=\emptyset$. Indeed, since $a^*\not\in\sigma$, Proposition \ref{Fed1} implies that there exists  $a_1\in A$ such that $a^*\ll_{C_{\rho}}a_1^*$ and $a_1^*\not\in\sigma$. Then $a_1\ll_{C_{\rho}}a$ and $a_1\in\sigma$. Since $a\in\mathbb{B}$, we get that $a_1\in\mathbb{B}$. Hence $a_1\ll_{\rho}a$ and $a_1\in\mathbb{B}\cap\sigma$. Then, by (\ref{*}), $\varphi(a)\cap\sigma'\not=\emptyset$. So, $\sigma'\in\iota_B(\varphi(a))$. Thus, $f_{\varphi}^{-1}(U)\subseteq\iota_B(\varphi(a))=V$. Note that, by Theorem \ref{iota}, $V$ is an open subset of $\Psi^a(B,\rho',\mathbb{B}')$.

Conversely, let $\sigma'\in\iota_B(\varphi(a))$ and $\s=f_{\varphi}(\sigma')$. Then $\varphi(a)\cap\sigma'\not=\emptyset$. We will prove that $a^*\not\in\sigma$. Suppose first that for every $e\ll a$, $\varphi(e)\cap\sigma'=\emptyset$. We have, by (M3), that $\varphi(a)=\bigvee\{\varphi(e)\st e\ll a\}$. Also, by Proposition \ref{clustideal}, $J_{\sigma\ap}=\mathbb{B}\ap\setminus\sigma'$ is a $\delta$-ideal. Since $\displaystyle\bigcup_{e\ll a}\varphi(e)\subseteq J_{\sigma\ap}$, we get that $\varphi(a)\subseteq J_{\sigma\ap}$, i.e. $\varphi(a)\cap\sigma'=\emptyset$, which is a contradiction. Hence, there  exists an $e\ll a$ such that $\varphi(e)\cap\sigma'\not=\emptyset$. Then $e\in\mathbb{B}$ (since $a\in\mathbb{B}$) and by (\ref{*}), $e\in\sigma\cap\mathbb{B}$. Since $e\ll_{\rho}a$, we have that $e(-\rho)a^*$. Using the fact that $e\in\mathbb{B}$, we get that $e(-C_{\rho})a^*$. Hence $a^*\not\in\sigma$. So, $\sigma\in \int(\lambda_A^g(a))=U$. Thus $\s\ap\in f_\p\inv(U)$. So, we have proved that

\begin{equation}\label{**}
f_{\varphi}^{-1}(\int(\lambda_A^g(a)))=\iota_B(\varphi(a)),\ \forall a\in\mathbb{B}.
\end{equation}

Now, using (\ref{eel}), we obtain that $f_{\varphi}$ is a continuous function. \sqs

\begin{pro}\label{functt}
For each $(A,\rho,\mathbb{B})\in|{\bf MDHLC}|$, set $\Delta^a(A,\rho,\mathbb{B})=\Psi^a(A,\rho,\mathbb{B})$ (see the text immediately after Theorem \ref{roeperl}
 for the notation $\Psi^a$), and for each $\MDHLC$-morphism $\varphi:(A,\rho,\mathbb{B})\lra (B,\rho',\mathbb{B}')$, put $\Delta^a(\varphi)=f_{\varphi}$ (see Proposition \ref{Delta^a} for the notation $f_\p$).
Then $\DE^a:\MDHLC\lra\HLC$ is a contravariant functor.
\end{pro}

\doc Let $(A,\rho,\mathbb{B})$ be a CLCA, $X=\DE^a(A,\rho,\mathbb{B})$ and $f=\DE^a(i_A)$. We will show that $f=id_X$. Indeed, by (\ref{Delata^ae}),
we have that for every $\s\in X$, $f(\s)\cap\BBBB=\{a\in\BBBB\st(\fa b\in A)[(a\ll_\rho b)\rightarrow(I_b\cap\s\nes)]\}$. By Proposition \ref{clustuniq}, it is enough to prove that
$f(\s)\cap\BBBB=\s\cap\BBBB$. Let $a\in f(\s)\cap\BBBB$. Suppose that $a\nin\s$. Then there exists $b\in\s$ such that $a(-C_\rho)b$. Thus $a\ll_\rho b^*$ and we get that
$I_{b^*}\cap\s\nes$. Let $c\in I_{b^*}\cap\s$. Then $c\in\BBBB$ and $c\ll_\rho b^*$. This implies that $c(-C_\rho)b$. Since $b,c\in\s$, we get a contradiction. So,
$f(\s)\cap\BBBB\sbe\s\cap\BBBB$. Conversely, let $a\in\s\cap\BBBB$. Let $b\in A$ and $a\ll_\rho b$. Then $a\in I_b\cap\s$, i.e. $I_b\cap\s\nes$. Thus $a\in f(\s)\cap\BBBB$. Hence $f(\s)\cap\BBBB=\s\cap\BBBB$. So, we have proved that $\DE^a(i_A)=id_X$.

Let now $\p_i\in\MDHLC((A_i,\rho_i,\BBBB_i),(A_{i+1},\rho_{i+1},\BBBB_{i+1}))$, where $i=1,2$, and $\p=\p_2\di\p_1$.
Set $f_i=\DE^a(\p_i)$, for $i=1,2$,   and let $f=\DE^a(\p)$. We will show that $f=f_1\circ f_2$. For $i=1,2,3$, set $X_i=\DE^a(A_i,\rho_i,\BBBB_i)$ and
$\ll_i=\ll_{\rho_i}$.
Let $\s_3\in X_3$ and set $\s_1\ap=f(\s_3)$. We have that $\s_1\ap\cap\BBBB_1=\{a_1\in\BBBB_1\st (\fa b_1\in A_1)[(a_1\ll_1 b_1)\rightarrow (\s_3\cap\bv\{\p_2(b_2)\st b_2\in\p_1(b_1)\}\nes)]\}=\{a_1\in\BBBB_1\st(\fa b_1\in A_1)[(a_1\ll_1 b_1)\rightarrow (\ex k\in\NNNN^+$ and $\ex c_1,\ldots, c_k\in\p_1(b_1)$ and $ \ex d_i\in\p_2(c_i),$
where $i=1,\ldots, k$, such that $\bv\{d_i\st i=1,\ldots, k\}\in\s_3)]\}=\{a_1\in\BBBB_1\st(\fa b_1\in A_1)[(a_1\ll_1 b_1)\rightarrow (\ex c\in\p_1(b_1)$ such that
$\p_2(c)\cap\s_3\nes)]\}=R$. Further, set $\s_2\ap=f_2(\s_3)$. Then we have that  $\s_2\ap\cap\BBBB_2=\{a_2\in\BBBB_2\st(\fa b_2\in A_2)[(a_2\ll_2 b_2)\rightarrow
(\ex c_2\in\p_2(b_2)\cap\s_3)]\}$. Now, $f_1(\s_2\ap)\cap\BBBB_1=\{a_1\in\BBBB_1\st(\fa b_1\in A_1)[(a_1\ll_1 b_1)\rightarrow(\ex c_2\in\p_1(b_1)\cap\s_2\ap)]\}=
\{a_1\in\BBBB_1\st(\fa b_1\in A_1)[(a_1\ll_1 b_1)\rightarrow(\ex c_2\in\p_1(b_1)$ such that $(\fa d_2\in A_2)((c_2\ll_2 d_2)\rightarrow(\p_2(d_2)\cap\s_3\nes)))]\}=R_{1,2}$.
By Proposition \ref{clustuniq}, it is enough to show that $R=R_{1,2}$. Let $a_1\in R$, $b_1\in A_1$ and $a_1\ll_1 b_1$. Then there exists $c_2\in\p_1(b_1)$ such that $\p_2(c_2)\cap\s_3\nes$. Let $d_2\in A_2$ and $c_2\ll_2 d_2$. Then $\p_2(d_2)\cap\s_3\nes$. Indeed, this follows from the facts that $\p_2(c_2)\sbe\p_2(d_2)$ and $\p_2(c_2)\cap\s_3\nes$. So, $a_1\in R_{1,2}$.
Conversely, let $a_1\in R_{1,2}$, $b_1\in A_1$ and $a_1\ll_1 b_1$. Then there exists $c_2\in\p_1(b_1)$ such that $(\fa d_2\in A_2)[(c_2\ll_2 d_2)\rightarrow
(\p_2(d_2)\cap\s_3\nes)]$. Since $\p_1(b_1)$ is a $\d$-ideal, there exists $c_2\ap\in\p_1(b_1)$ such that $c_2\ll_2 c_2\ap$. Then $\p_2(c_2\ap)\cap\s_3\nes$. Therefore,
$a_1\in R$. So, we have proved that $f=f_1\circ f_2$. All this shows that $\DE^a$ is a contravariant functor. \sqs

\begin{pro}\label{widetilde}
If $\varphi:(A,\rho,\BBBB)\longrightarrow (B,\eta,\BBBB\ap)$ is an {\bf LCA}-isomorphism then the multi-valued map
 $\widetilde{\varphi}:(A,\rho,\BBBB)\longrightarrow (B,\eta,\BBBB\ap)$, where $\widetilde{\varphi}(a)=I_{\varphi(a)}$, is a {\bf MDHLC}-isomorphism.
\end{pro}

\doc It is obvious that $\widetilde{\varphi}$ satisfies conditions (M1) and (M4). Further, we have that $\widetilde{\varphi}(a\wedge b)=I_{\varphi(a\wedge b)}=I_{\varphi(a)\wedge\varphi(b)}=I_{\varphi(a)}\cap I_{\varphi(b)}=\widetilde{\varphi}(a)\we\widetilde{\varphi}(b)$. So, condition (M2) is fulfilled.

We will prove  that for every $a\in A$, $\widetilde{\varphi}(a)=\bigvee\{\widetilde{\varphi}(b)\st b\in\mathbb{B},\ b\ll a\}$, i.e. $$I_{\varphi(a)}=\bigvee\{I_{\varphi(b)}\st b\in\mathbb{B},\ b\ll a\}.$$ Indeed, let $b\in\mathbb{B}$ and $b\ll a$. Then $\varphi(b)\ll\varphi(a)$. Hence $I_{\varphi(b)}\subseteq I_{\varphi(a)}$. Therefore $\bigvee\{I_{\varphi(b)}\st b\in\mathbb{B},\ b\ll a\}\subseteq I_{\varphi(a)}$. Conversely, let $c'\in I_{\varphi(a)}$. Then $c'\in\mathbb{B}'$ and $c'\ll\varphi(a)$. There exists  $c''\in\mathbb{B}'$ such that $c'\ll c''\ll\varphi(a)$. There exists  $c\in\mathbb{B}$ such that $c''=\varphi(c)$. Then $\varphi(c)\ll\varphi(a)$; hence $c\ll a$ and $\varphi(c)=c''\gg c'$. Therefore $c'\in I_{\varphi(c)}$, where $c\in\mathbb{B}$ and $c\ll a$. Thus, $I_{\varphi(a)}\sbe\bigcup\{I_{\varphi(b)}\st b\in\mathbb{B},\ b\ll a\}\sbe\bigvee\{I_{\varphi(b)}\st b\in\mathbb{B},\ b\ll a\}$.
So, condition (M3) is also fulfilled.

We will now verify (M5). Let $a_i,b_i\in\mathbb{B}$ and $a_i\ll b_i$, where $i=1,2$. We will prove that $\widetilde{\varphi}(a_1\vee a_2)\subseteq\widetilde{\varphi}(b_1)\vee\widetilde{\varphi}(b_2)$, i.e. $I_{\varphi(a_1\vee a_2)}\subseteq I_{\varphi(b_1)}\vee I_{\varphi(b_2)}$. Indeed, let $c\in\mathbb{B}$ and $c\ll\varphi(a_1\vee a_2)$. Then $c\ll\varphi(a_1)\vee\varphi(a_2)$. We have that
 $c\wedge\varphi(a_1)\leq\varphi(a_1)\ll\varphi(b_1)$, $c\wedge\varphi(a_2)\leq\varphi(a_2)\ll\varphi(b_2)$ and $c=(c\wedge\varphi(a_1))\vee(c\wedge\varphi(a_2))$. Set $d_i=c\wedge\varphi(a_i)$, for $i=1,2$. Then $d_i\ll\varphi(b_i)$, i.e. $d_i\in I_{\varphi(b_i)}$, for $i=1,2$, and $c=d_1\vee d_2$. Hence $c\in I_{\varphi(b_1)}\vee I_{\varphi(b_2)}$. So, $I_{\varphi(a_1\vee a_2)}\subseteq I_{\varphi(b_1)}\vee I_{\varphi(b_2)}$, i.e. $\widetilde{\varphi}(a_1\vee a_2)\subseteq\widetilde{\varphi}(b_1)\vee\widetilde{\varphi}(b_2)$.

We will show that condition (M6) is satisfied, i.e. that $\bigcup\{\widetilde{\varphi}(a)\st a\in\mathbb{B}\}=\mathbb{B}'$ holds. Indeed, let $b'\in\mathbb{B}'$. Then there exists  $b''\in\mathbb{B}'$ such that $b'\ll b''$. There exists an $a\in\mathbb{B}$ such that $b''=\varphi(a)$. Then $b'\in I_{\varphi(a)}=\widetilde{\varphi}(a)$.

Hence, $\widetilde{\varphi}$ is an {\bf MDHLC}-morphism. Analogously, we obtain that $\widetilde{\varphi^{-1}}$  is an {\bf MDHLC}-morphism.

We will prove that $\widetilde{\varphi}\diamond\widetilde{\varphi^{-1}}=i_B$ and $\widetilde{\varphi^{-1}}\diamond\widetilde{\varphi}=i_A$. Indeed, $(\widetilde{\varphi^{-1}}\diamond\widetilde{\varphi})(a)= \bigvee\{\widetilde{\varphi^{-1}}(b)\st b\in\widetilde{\varphi}(a)\}= \bigvee\{I_{\varphi^{-1}(b)}\st b\in I_{\varphi(a)}\}$ and $i_A(a)=I_a$ for every $a\in A$. So, we have to prove that $I_a=\bigvee\{I_{\varphi^{-1}(b)}\st b\in I_{\varphi(a)}\}$. Indeed, let $c\in I_a$. Then $c\in\mathbb{B}$ and $c\ll a$. Hence there exists  $d\in\mathbb{B}$ such that $c\ll d\ll a$. Set $b=\varphi(d)$. Then $b\ll\varphi(a)$, i.e. $b\in I_{\varphi(a)}$. Also, $c\ll d=\varphi^{-1}(\varphi(d))=\varphi^{-1}(b)$, i.e. $c\in I_{\varphi^{-1}(b)}$. Hence $I_a\subseteq\bigcup\{I_{\varphi^{-1}(b)}\st b\in I_{\varphi(a)}\}$. Conversely, let $c\in I_{\varphi^{-1}(b)}$, where $b\in I_{\varphi(a)}$. Then $c\ll\varphi^{-1}(b)$ and $b\ll\varphi(a)$. Since $\varphi^{-1}(b)\ll\varphi^{-1}\varphi(a)=a$, we get that $c\ll a$, i.e. $c\in I_a$. So, $\bigcup\{I_{\varphi^{-1}(b)}\st b\in I_{\varphi(a)}\}\subseteq I_a$. Hence $I_a=\bigcup\{I_{\varphi^{-1}(b)}\st b\in I_{\varphi(a)}\}$. Then $I_a=\bigvee\{I_{\varphi^{-1}(b)}\st b\in I_{\varphi(a)}\}$. So, $\widetilde{\varphi^{-1}}\diamond\widetilde{\varphi}=i_A$. Analogously, we get that $\widetilde{\varphi}\diamond\widetilde{\varphi^{-1}}=i_B$.

Therefore, $\widetilde{\varphi}$ is a {\bf MDHLC}-isomorphism. \sqs

\begin{pro}\label{diagrama1}
The identity functor $Id_{\bf MDHLC}$ and the functor $\Delta^t\circ\Delta^a$
are naturally isomorphic.
\end{pro}

\doc Let $\varphi\in{\bf MDHLC}((A,\rho,\mathbb{B}),(B,\eta,\mathbb{B}\ap))$. We have to show that $\widetilde{\lambda_B^g}\diamond\varphi=\Delta^t(\Delta^a(\varphi))\diamond\widetilde{\lambda_A^g}$, where $\widetilde{\lambda_A^g}(a)=I_{\lambda_A^g(a)}$ (see \ref{widetilde}).
(Note that, by  (\ref{hapisomn}), $\lag$ and $\lbg$ are LCA-isomorphisms and, hence, by Proposition \ref{widetilde}, $\wt{\lag}$ and $\wt{\lbg}$ are $\MDHLC$-isomorphisms.)

Set $\Delta^a(A,\rho,\BBBB)=X$, $\Delta^a(B,\eta,\BBBB\ap)=Y$ and
 $\p\ap=\Delta^t(\Delta^a(\varphi))$ $(=\Delta^t(f_{\varphi}))$ Hence $\varphi\ap:(RC(X),\rho_X,CR(X))\longrightarrow (RC(Y),\rho_Y,CR(Y))$.
Then, for each $F\in RC(X)$, $\varphi\ap(F)=\{G\in CR(Y)\st G\subseteq f_{\varphi}^{-1}(\int (F))\}$. Hence, for every $a\in A$,

\noindent$\begin{array}{rl}
(\varphi\ap\diamond\widetilde{\lag})(a)&=\bigvee\{\varphi\ap(b)\st b\in\widetilde{\lag}(a)\}=\bigvee\{\varphi\ap(b)\st b\in I_{\lambda_A^g(a)}\}\\
&=\bigvee\{\varphi\ap(G)\st G\in CR(X),\ G\ll\lambda_A^g(a)\}\\
&=\bigvee\{\varphi\ap(G)\st G\in CR(X),\ G\subseteq \int(\lambda_A^g(a))\}\\
&=\bigvee\{\{H\in CR(Y)\st H\subseteq f_{\varphi}^{-1}(\int G)\}\st G\in CR(X),\ G\subseteq \int\lambda_A^g(a)\}\\
&=\bigvee\{\{\lambda_B^g(b\ap)\st b\ap\in\mathbb{B}\ap,\ \lambda_B^g(b\ap)\subseteq f_{\varphi}^{-1}(\int(\lambda_A^g(c)))\}\st c\in\mathbb{B},\ c\ll_{\rho}a\}.
\end{array}$

\noindent Since, by (\ref{**}), $f_{\varphi}^{-1}(\int(\lambda_A^g(a)))=\iota_B(\varphi(a))$, we get that

\noindent$\begin{array}{rl}
(\varphi\ap\diamond\widetilde{\lag})(a)&=\bigvee\{\{\lambda_B^g(b\ap)\st b\ap\in\mathbb{B}\ap,\ \lambda_B^g(b\ap)\subseteq\iota_B(\varphi(c))\}\st c\in\mathbb{B},\ c\ll_{\rho}a\}\\
&=\bigvee\{\{\lambda_B^g(b\ap)\st b\ap\in\varphi(c)\}\st c\in\mathbb{B},\ c\ll a\}.
\end{array}$

\noindent The last equality follows from the fact that for every $b\ap\in\BBBB\ap$, $\lbg(b\ap)$ is compact and hence there exist
$b_1\ap,\ldots, b_n\ap\in\p(c)$ such that $\lbg(b\ap)\le\bv\{\lbg(b_i\ap)\st i=1,\ldots,n\}$; conversely, for every $b\ap\in\p(c)$, $\lbg(b\ap)\sbe\iota_B(\p(c))$.

Further,
$$\begin{array}{rl}
(\widetilde{\lbg}\diamond\varphi)(a)&=\bigvee\{\widetilde{\lbg}(b)\st b\in\varphi(a)\}\\
&= \bigvee\{I_{\lambda_B^g(b)}\st b\in\varphi(a)\}\\
&= \bigvee\{\{\lambda_B^g(b\ap)\st b\ap\in\mathbb{B},\ b\ap\ll b\}\st b\in\varphi(a)\}.
\end{array}$$

 Hence
 $$\begin{array}{rl}
 (\varphi\ap\diamond\widetilde{\lag})(a)&=\{\lambda_B^g(b_1\ap\vee\ldots\vee b_k\ap)\st b_i\ap\in\varphi(c_i), c_i\in\mathbb{B}, c_i\ll a, k\in\NNNN^+, i=1,\ldots,k \}\\
  &= \{\lbg(b\ap)\st b\ap\in\p(c), c\in\BBBB, c\ll a\}.
  \end{array}$$
  and
  $$\begin{array}{rl}
  (\widetilde{\lbg}\diamond\varphi)(a)&=\{\lambda_B^g(b_1\ap\vee\cdots\vee b_k\ap)\st b_i\ap\ll b_i,\ b_i\in\varphi(a),\ k\in\NNNN^+,\ i=1,\ldots,k \}\\
  &= \{\lbg(b\ap)\st b\ap\in\p(a)\}=\lbg(\p(a)).
  \end{array}$$

Let $b\ap$ be such that $\lambda_B^g(b\ap)\in(\varphi\ap\diamond\widetilde{\lag})(a)$, i.e.  $b\ap\in\varphi(c)$, where $c\in\mathbb{B}$, $c\ll a$. Since $\varphi(c)\subseteq\varphi(a)$, we get that $\lambda_B^g(b\ap)\in(\widetilde{\lbg}\diamond\varphi)(a)$.

Conversely, let $b\ap$ be such that $\lambda_B^g(b\ap)\in(\widetilde{\lbg}\diamond\varphi)(a)$, i.e.  $b\ap\in\varphi(a)$. By (M3), $\varphi(a)=\bigvee\{\varphi(c)\st c\in\mathbb{B},\ c\ll a\}= \{d_1\vee\cdots\vee d_k\st k\in\mathbb{N}^+,\ d_i\in\varphi(c_i),\ c_i\ll a,\ c_i\in\mathbb{B}\}$. Hence $b\ap=d_1\vee\cdots\vee d_k$, $d_i\in\varphi(c_i)$, $c_i\in\mathbb{B}$, $c_i\ll a$, for every $i=1,\ldots,k$. Set $c=\bv\{c_i\st i=1,\ldots,k\}$.
Then $c\ll a$, $c\in\BBBB$ and $d_i\in\p(c)$  for every $i=1,\ldots,k$.
Hence $b\ap\in\p(c)$.
Thus, $\lambda_B^g(b\ap)\in(\varphi\ap\diamond\widetilde{\lag})(a)$.

Hence, $\widetilde{\lbg}\diamond\varphi=\Delta^t(\Delta^a(\varphi))\diamond\widetilde{\lag}$. \sqs

\begin{pro}\label{diagrama2}
The identity functor $Id_{\bf HLC}$ and the functor $\Delta^a\circ\Delta^t$
are naturally isomorphic.
\end{pro}

\doc Let $f\in{\bf HLC}(X,Y)$. We have to show that
$t_Y\circ f=\Delta^a(\Delta^t(f))\circ t_X$, where $t_X(x)=\sigma_x$ for every $x\in X$. (Recall that, by (\ref{homeo}), $t_X$ and $t_Y$ are homeomorphisms.)
Set $f'=\Delta^a(\Delta^t(f))$ $(=\Delta^a(\p_f))$. Then, for each $\sigma\in\Delta^a(\Delta^t(X))$, we have that $f'(\sigma)=\sigma'$,
where $f'(\sigma)\cap CR(Y)=\{G\in CR(Y)\st(\forall H\in RC(Y))((G\subseteq \int
(H))\rightarrow(\varphi_f(H)\cap\sigma\not=\emptyset))\}$.

Now, for every $x\in X$, $(f'\circ t_X)(x)=f'(\sigma_x)=\sigma'$, where
$\sigma'\cap CR(Y)=\{G\in CR(Y)\st$ $(\forall H\in RC(Y))((G\subseteq \int (H))\rightarrow(\exists F\in RC(X)\mbox{ such that } x\in F\mbox{ and } F\in\varphi_f(H)))\}.$ Hence $\sigma'\cap CR(Y)=\{G\in CR(Y)\st (\forall H\in RC(Y))((G\subseteq \int (H))\rightarrow(\exists F\in RC(X)\mbox{ such that } x\in F\subseteq f^{-1}(\int (H))))\}$.

Further, $(t_Y\circ f)(x)=\sigma_{f(x)}$, where $\sigma_{f(x)}\cap CR(Y)=\{G\in CR(Y)\st f(x)\in G\}$.

Let $G\in\sigma_{f(x)}\cap CR(Y)$. Then $f(x)\in G$. We will prove that $G\in\sigma'$. Let $H\in CR(Y)$ and $G\subseteq \int (H)$.
We will prove that there exists an $F\in RC(X)$ such that $x\in F\subseteq f^{-1}(\int (H))$.
Indeed, $f(x)\in G\subseteq \int (H)$. Since $f$ is continuous, there exists an open $U\subseteq X$ such that $x\in U$ and $f(U)\subseteq \int (H)$.
Since $X$ is a locally compact $T_2$-space, there exists an $F\in CR(X)$ such that $x\in F\subseteq U$.
Then $f(F)\subseteq f(U)\subseteq \int (H)$, i.e. $F\subseteq f^{-1}(\int (H))$. So, $G\in\sigma'\cap CR(Y)$.
Hence $\sigma_{f(x)}\cap CR(Y)\subseteq\sigma'\cap CR(Y)$.

Conversely, let $G\in CR(Y)\cap\sigma'$. We will prove that $f(x)\in G$. Indeed, suppose that $f(x)\not\in G$. Then there exists an $H\in CR(Y)$ such that $G\subseteq \int (H)\subseteq Y\setminus\{f(x)\}$. We have that there exists an $F\in CR(X)$ such that $x\in F\subseteq f^{-1}(\int (H))$. Then $f(x)\in \int (H)$, which is a contradiction. So, $f(x)\in G$. Hence $\sigma_{f(x)}\cap CR(Y)\supseteq\sigma'\cap CR(Y)$.

We get that $\sigma_{f(x)}\cap CR(Y)=\sigma'\cap CR(Y)$. Then, by Proposition \ref{clustuniq}, $\sigma_{f(x)}\equiv\sigma'$. So, $t_Y\circ f=\Delta^a(\Delta^t(f))\circ t_X$. \sqs

The next theorem, which is the main result of this paper, follows from Theorem \ref{roeperl} and Propositions \ref{functa},  \ref{functt}, \ref{diagrama1},  \ref{diagrama2}.

\begin{theorem}\label{main}{\rm (The Main Theorem)}
The categories {\bf HLC} and {\bf MDHLC} are dually equivalent.
\end{theorem}

\end{document}